\documentclass[12pt]{amsart}
\begingroup
\swapnumbers
\theoremstyle{plain}

\newtheorem{thm}{Theorem}[section]
\newtheorem{cor}[thm]{Corollary}

\newtheorem{lem}[thm]{Lemma}

\newtheorem{prob}[thm]{Problem}

\theoremstyle{definition}
\newtheorem{dfn}[thm]{Definition}
\newtheorem{rmk}[thm]{Remark}
\newtheorem{rmks}[thm]{Remarks}

\endgroup

\numberwithin{equation}{section}

\newcommand {\Bh}{B(H)}  
\newcommand {\Bx}{B(X)}  
\newcommand {\By}{B(Y)}  
\newcommand {\N}{\mathbb{N}} 
\newcommand {\Z}{\mathbb{Z}}        
\newcommand {\C}{\mathbb{C}} 
\newcommand {\D}{\mathbb{D}} 
\newcommand {\T}{\mathbb{T}} 

\newcommand {\sm}{\sigma}
\def\e{\varepsilon}
\newcommand {\ld}{\lambda} 
\def\AA{{\mathcal A}}
\def\BB{{\mathcal B}}
\def\CCC{{\mathcal C}}
\newcommand {\et}{\mathcal{E}(\T)}
\def\al{\alpha}

\newcommand {\be}{\begin{equation}}
\newcommand {\ee}{\end{equation}}
\newcommand {\beq}{\begin{eqnarray*}}
\newcommand {\eeq}{\end{eqnarray*}}

\begin{document}

\title
{Growth conditions and inverse producing extensions}

\author{C. Badea}
\address{D\'epartement de Math\'ematiques, UMR CNRS no. 8524, 
Universit\'e Lille I,
  F--59655 Villeneuve d'Ascq, France}
\email{Catalin.Badea@math.univ-lille1.fr}
\thanks{C.B. was supported by
the EU Research Training Network 'Analysis and Operators' 
HPRN-CT-2000-00116, the EU programme 
EUROMMAT 
ICA1-CT-2000-70022 and by CNRS (France) through a 'd\'el\'egation CNRS' 
at Institut Henri Poincar\'e, Paris} 

\author{ V. M{\"u}ller}
\address{Institut of Mathematics AV CR, Zitna 25, 115 67 Prague 1, Czech
                   Republic}
\email{muller@math.cas.cz}
\thanks{V.M. was supported by grant No. 201/03/0041 of GA \v CR}

\subjclass[2000]{47A20, 47A05, 46J05, 47B40}
\keywords{invertible extensions, growth conditions, subscalar operators}

\date{}

\begin{abstract}
We study the invertibility of Banach algebras elements in their
extensions, and invertible extensions 
of Banach and Hilbert space 
operators with prescribed growth conditions for the norm of inverses. 
As applications, the solutions of two open problems are obtained. 
In the first one we give a characterization of $\et$-subscalar 
operators in terms of growth conditions. 
In the second one we show that operators satisfying a Beurling-type growth
condition possess Bishop's property $(\beta)$. 
Other applications are also given.
\end{abstract}

\maketitle


\section{Introduction}
\subsection{Preamble.} A bounded linear operator 
can be made 'nicer' by an extension or a 
dilation to a
larger space. 
One example, \cite{szfo}, is
the celebrated Sz.-Nagy dilation theorem (every Hilbert space contraction has a
unitary dilation), or its extension variant (every Hilbert space contraction has a
coisometric extension). A Banach space example is a result 
due to R.G.~Douglas, \cite{douglas},
stating that a Banach space isometry has an extension to a
surjective isometry. Douglas' construction is hilbertian, in the
sense that if the given operator acts on a Hilbert space, then its extension, a
unitary operator, acts also on a Hilbert space. 
In the framework of Banach algebras, 
a classical result of R.F.~Arens, \cite{arens}, states that if an element $u$ 
of a commutative unital Banach algebra $\AA$ is not a 
topological divisor of zero, then $u$ is 
invertible in a commutative unital Banach algebra 
containing $\AA$. Other such examples, related to the topic of the present
paper, can be found in
\cite{stroescu,read1,muller,read2,read3,bp,batty/yeates}.
\subsection{Motivation.} The aim of this paper is to study the
invertibility of Banach algebras elements in their extensions,
and invertible 
extensions of Banach or Hilbert space 
operators with prescribed growth conditions for the norm of inverses. We obtain,
among other things, generalizations of the above mentioned results of
Douglas and Arens.

Our investigations were also motivated by two open problems, which 
will be solved positively in this paper. The first one is 
due to K.B.~Laursen and
M.M.~Neumann \cite[Problem 6.1.15]{laursen/neumann} 
and M.~Didas \cite{didas} and asks for a 
characterization in terms of growth conditions 
of $\et$-subscalar operators, i.e., of operators which are similar to 
restrictions of 
$\et$-scalar operators to closed invariant subspaces. 

The second open problem asks \cite{millers/neumann} if
operators $T\in B(X)$ satisfying the Beurling-type condition
\begin{equation}\label{eq:beurling}
\sum_{n=1}^{\infty}\frac{\log \max(\|T^n\|,m(T^n)^{-1})}{n^2} < \infty 
\end{equation}
possess Bishop's property $(\beta)$ ; see (\ref{eq:1.2}) for the definition of
the minimum modulus $m(T^n)$ and Section \ref{sect:applications} for 
the definition of property $(\beta)$. 
\subsection{Organization of the paper.}
Our first result in the second 
section is a refinement of the Arens construction. We consider
the invertibility of an element $u$ of a Banach algebra $\AA$ in an
extension of $\AA$ with prescribed growth
conditions for $\|u^{-k}\|$, $k \ge 1$. We then 
consider extensions of Banach space operators. 
We use a method due to one of the authors 
\cite{muller} to pass from
the Banach algebra case to the case of $B(X)$. 

In Section 3
we use an idea of Batty and Yeates \cite{batty/yeates} 
to show that, given a real
number $p\ge 1$ and $T \in \Bx$, there is an isomorphic embedding $\pi : X
\mapsto Y$ and an invertible operator 
$S \in \By$ with prescribed growth conditions for $\|S^{-k}\|$, $k \ge 1$, 
such that $T$ is similar
to the restriction of $S$ to $\pi(X)$. Moreover, 
the space $Y$ may be obtained from
$X$ as a quotient of a subspace of an ultraproduct of spaces of the form
$L_p(X)$ (i.e., a $SQ_p(X)$-space). 
In particular, if $p=2$ and $X$ is a Hilbert space, then so is $Y$. 

In the last section we consider several applications. A characterization for 
$\et$-subscalar operators is 
given in Theorem \ref{subscalar}. 
The question from \cite{millers/neumann} concerning
operators satisfying the Beurling-type 
condition (\ref{eq:beurling}) is positively
answered in Theorem \ref{beta}. We then consider operators satisfying some 
exponential growth conditions. Other applications 
concerning operators with countable spectrum and Hilbert space contractions with
spectrum a Carleson set are given.
\subsection{Notation and terminology.} We recall now some known facts and 
introduce some notation. All other undefined terms are classical or 
will be defined in
Section~\ref{sect:applications}.

{\bf Banach algebras.} All Banach algebras are considered 
to be complex and with unit. Let $u$ be an element 
of a Banach algebra
$\AA$. We write 
$$
d^{\AA}(u) = \inf
\{\|ux\|: x\in \AA, \|x\| = 1\}.
$$ 
If no confusion can arise then we omit the upper index and 
write simply $d(u)$ instead of
$d^{\AA}(u)$.  

Let $\AA,\BB$ be commutative Banach algebras. We say that $\BB$
is an extension of $\AA$ if there exists an isometrical unit
preserving homomorphism
$\rho:\AA\to\BB$.
If we identify $\AA$ with the image $\rho(\AA)$ we can consider
$\AA$ as a closed subalgebra of $\BB$ and write simply $\AA\subset\BB$.

{\bf Operators.} In this paper $X$ (and $Y$) will denote complex
Banach spaces and $H$ (and $K$) will 
denote Hilbert spaces. Denote by $B(X)$ the algebra of all bounded
linear operators on the Banach space $X$.
By an operator we always mean a
bounded linear operator. Note that for an
operator $T\in B(X)$ we can express the quantity $d^{B(X)}(T)$
in a more convenient way by 
\begin{equation}\label{eq:1.2}
m(T) := d^{B(X)}(T) = \inf\{\|Tx\|:x\in X,
\|x\|=1\}.
\end{equation}
This quantity is called the \emph{minimum modulus} of $T$ (\cite{gindler}) or 
the \emph{lower bound} of $T$ (\cite{laursen/neumann}). 

We denote by $\sm(T)$ and $\sm_{ap}(T)$ the spectrum 
and the approximate
point spectrum of a bounded linear operator $T\in\Bx$,
respectively. The latter is given by
$$\sm_{ap}(T) = \bigl\{\ld \in \C : \inf\{\|(T-\ld)x\| : \|x\|
=1\} = 0\bigr\}. $$
Note that $m(T) > 0$ if and only if 
$T\in B(X)$ is one-to-one and of closed range. 
If T is a Hilbert space operator, then $\sm_{ap}(T)$ coincides
with the left spectrum and 
$m(T) > 0$ if and only if $T$ is left invertible. 

We say that $S \in \By$ is an extension of $T\in \Bx$ if  there is an isometry
$\pi : X \to Y$ such that $S\pi = \pi T$. 
We can also consider $X$ as a subspace of $Y$ and write $T=S_{\mid X}$.

{\bf Banach spaces of class $SQ_p$.}
Let $p\geq 1$ be a real number. A Banach space $E$ is said to be a 
$SQ_p$-space if it is a quotient of a subspace of an $L_p$-space. 

Let $X$ be a Banach space. A Banach space $E$ is said to be a 
$SQ_p(X)$-space if it is (isometric to) a quotient of a subspace of 
an ultraproduct of spaces of the form $L_p(\Omega,\mu,X)$, for some 
measure spaces $(\Omega,\mu)$. 
Since ultraproducts of $L_p$-spaces are $L_p$-spaces, the latter
definition is consistent with the former one. 
Note that any Banach space is isometric to a
subspace (resp. a quotient) of an $L_{\infty}$-space (resp. an $L_1$-space). 
Also, if $H$ is a
Hilbert space, then each $SQ_2(H)$-space is a Hilbert space too.

$SQ_p(X)$-spaces are characterized 
by a theorem of R.~Hernandez \cite{hernandez} (for $X=\C$ this goes back to 
\cite{kwapien}). 
See also \cite{pisier:indiana} (and \cite[Theorem 3.2]{lemerdy}) 
for a different proof using 
$p$-completely bounded maps. Namely, $E$ is a $SQ_p(X)$-space 
if and only if 
$$\|a\|_{p,E} \leq \|a\|_{p,X}$$
for each $n\geq 1$ and each matrix  $a = [a_{ij}]\in M_{n}(\C)$. Here 
$$\|[a_{ij}]\|_{p,Y} =  \sup 
\left[\left(\sum_i\Bigl\|\sum_ja_{ij}y_j\Bigr\|^p\right)^{1/p}\right] ,$$
where the supremum runs over all $n$-tuples $(y_1, \cdots , y_n)$ in 
$Y$ which satisfy $\sum\|y_j\|^p \leq 1$.

{\bf Nearness.} 
Let $p \ge 1$ and $\beta : \N \to (0,\infty)$. 
Let $X$ be a subspace of $Y$. Two operators $T$ and $C$ in $\By$ 
are said to be $(\beta,p)$-\emph{near modulo} $X$ if 
for every $N \in \N$ and for all $x_1, \ldots , x_N \in X$ we have
\begin{equation}\label{eq:mod}
\left\|\sum_{n=1}^{N} (T^n - C^n)x_n\right\| \leq  
\left(\sum_{n=1}^{N}\beta(n)^p\| x_n\|^p\right)^{1/p} .
\end{equation}
For a constant weight function $\beta(n) \equiv s$ and for $p =2$ this definition
was introduced and studied 
in \cite{badea,badea2} under the name of quadratic nearness.

Note that if
$p=1$, and if the operators $T,C\in\By$ verify
$\|T^{n} - C^n\|\le \beta(n)$ for all $n\ge 1$, then (\ref{eq:mod}) 
holds for every
$x_n \in Y$.

\section{A refinement of the Arens construction}
The result of R.F.~Arens \cite{arens} implies that if
$\AA$ is a commutative Banach algebra and $d^\AA(u)>0$, 
then there exists a commutative extension
$\BB\supset\AA$ such that $u$ is invertible in $\BB$. It follows
from the Arens construction that $\|u^{-k}\|\le
\bigl(d^\AA(u)\bigr)^{-k}\quad(k\ge 1)$.
The following theorem gives a necessary and sufficient condition for having
invertible extensions of Banach algebra elements 
with prescribed growth conditions for the norm of
inverses.
\begin{thm}\label{thm:arensCNS}
Let $u$ be an element of a commutative Banach
algebra $\AA$. Let $(c_j)_{j=1}^\infty$ be a sequence of positive numbers 
which is submultiplicative, i.e., $c_{i+j}\le
c_ic_j$ for all $i,j\ge 1$. Then there is a
commutative extension $\BB\supset\AA$ such that $u$ is
invertible in $\BB$ and $\|u^{-j}\|\le c_j\quad(j\ge 1)$ if and only if 
we have
$$
\|a_0\| \leq \sum_{j=1}^\infty c_j\|a_j-a_{j-1}u\|
$$
for every sequence $(a_j)_{j=0}^{\infty}$ in $\AA$ of finite support.
\end{thm}
\begin{proof}
Suppose that $\BB\supset\AA$ is a commutative extension with all the
required properties. Let $(a_j)_{j=0}^{\infty}$ be a sequence in $\AA$ such that
$a_j=0$ for $j \geq n$. Write $f_j = a_j - a_{j-1}u$. Then 
\beq
\|a_0\|_{\AA}  & = &  \|a_0\|_{\BB} = \|u^{-n}u^na_0\|\\
        & = & \left\|-u^{-n}\left(\sum_{j=1}^nu^{n-j}f_j\right)\right\|\\
        & = & \Bigl\|\sum_{j=1}^n u^{-j}f_j\Bigr\| 
\leq \sum_{j=1}^nc_j\|f_j\| .
\eeq

For the converse, set formally $c_0=1$.
Consider the algebra $\CCC$ of all power
series $\sum_{i=0}^\infty a_i x^i$ 
in one variable $x$ with coefficients $a_i\in\AA$ such that
$$
\Bigl\|\sum_{i=0}^\infty a_i x^i\Bigr\|=\sum_{i=0}^\infty \|a_i\|c_i
<\infty.
$$
With the multiplication given by
$$
\left(\sum_{i=0}^\infty a_i x^i\right)\cdot \left(\sum_{j=0}^\infty
a'_j x^j\right) = \sum_{k=0}^\infty x^k \left(\sum_{i+j=k} a_i
a'_j\right), 
$$
$\CCC$  is a commutative Banach algebra containing $\AA$ as
subalgebra of constants. Let $J$ be the closed ideal generated by the
element $1-ux$ and set $\BB=\CCC/J$. 
Let $\rho:\AA\to\BB$ be the composition of the embedding
$\AA\to\CCC$ and the canonical homomorphism
$\CCC\to\BB=\CCC/J$.
Then
$$
\rho(u)\cdot(x+J)=(u+J)(x+J)=1_\AA+J=1_\BB,
$$
and so $\rho(u)$ is invertible in $\BB$ with the inverse $x+J$.
We have $\|(x+J)^n\|_\BB\le \|x^n\|_\CCC=c_n$ for all $n\ge 1$.

It is sufficient to show that $\rho$ is an isometry,
i.e., that for each $a\in \AA$ we have
$\|a\|_\AA=\|\rho(a)\|_\BB.$

Obviously, $\|\rho(a)\|_\BB=\inf_{c\in\CCC} \|a+(1-ux)c\|
\le\|a\|_\AA$. 

Suppose on the contrary that there is an $a\in \AA$ such that
$\|\rho(a)\|_\BB<\|a\|_\AA$. Thus there are elements $a_j\in\AA$
such that
\beq
\|a\|_\AA&>&\Bigl\|a+(1-ux)\sum_{j=0}^\infty a_jx^j\Bigr\|_\CCC\\
&=&
\|a-a_0\|_\AA+\sum_{j=1}^\infty c_j\cdot \|a_j-a_{j-1}u\|_\AA\\
&\ge&
\|a\|-\|a_0\|+\sum_{j=1}^\infty c_j\cdot \|a_j-a_{j-1}u\|.
\eeq
Thus $\|a_0\|>\sum_{j=1}^\infty c_j\|f_j\|$, where
$f_j=a_j-a_{j-1}u$. Moreover,
we may assume that only a finite number of elements $a_j$ are
non-zero. This contradicts to our assumption.
\end{proof}
We introduce the following 
definition.
\begin{dfn}Let $u$ be an element of a Banach algebra $\AA$. Let
$(c_j)_{j=1}^\infty$ be a sequence of positive numbers 
which is submultiplicative, i.e., $c_{i+j}\le
c_ic_j$ for all $i,j\ge 1$. We say
that $(c_j)$ \emph{satisfies condition} $(*)$ \emph{for} $u\in \AA$ if there
exists an increasing sequence $(k_n)$ of integers such that 
$0=k_0<k_1< k_2<\cdots$ and
$$
(*) \quad c_j\ge \Bigl(d(u^{k_1}) d(u^{k_2-k_1})\cdots
d(u^{k_{n+1}-k_n})\Bigr)^{-1} \|u^{k_{n+1}-j}\|
$$
for all $n\ge 0$ and $j$ satisfying $k_n<j\le k_{n+1}$.
\end{dfn}
\begin{thm}\label{thm:arens}
Let $u$ be an element of a commutative Banach
algebra $\AA$. Let $(c_j)$ be a sequence of positive numbers
satisfying condition $(*)$ for $u\in\AA$. Then there is a
commutative extension $\BB\supset\AA$ such that $u$ is
invertible in $\BB$ and $\|u^{-j}\|\le c_j\quad(j\ge 1)$.
\end{thm}
\begin{proof}
Set formally $c_0=1$.
Let $(a_j)_{j=0}^{\infty}$ be a sequence in $\AA$ of finite support. 
Write $f_j = a_j - a_{j-1}u$. 

We verify the condition of Theorem \ref{thm:arensCNS}.
We have
\beq
\|a_0\|&\le&
d(u^{k_1})^{-1}\|a_0u^{k_1}\|\\
&\le&
d(u^{k_1})^{-1}\Bigl(\|a_0u^{k_1}-a_1u^{k_1-1}\|+\cdots+\|a_{k_1-1}u-a_{k_1}\|+
\|a_{k_1}\|\Bigr)\\
&\le&
d(u^{k_1})^{-1}\Bigl(\|f_1\|\cdot\|u^{k_1-1}\|+\|f_2\|\cdot\|u^{k_1-2}\|+
\cdots+\|f_{k_1}\|
\Bigr)\\
& \hskip0.5cm&+d(u^{k_1})^{-1}d(u^{k_2-k_1})^{-1}\|a_{k_1}u^{k_2-k_1}\|\\
&\le&
\sum_{j=1}^{k_1}c_j\|f_j\|+d(u^{k_1})^{-1}d(u^{k_2-k_1})^{-1}
\Bigl(\bigl\|a_{k_1}u^{k_2-k_1}-a_{k_1+1}u^{k_2-k_1-1}\bigr\|\\
&\hskip0.5cm &+\cdots +
\|a_{k_2-1}u-a_{k_2}\|+\|a_{k_2}\|\Bigr)\\
&\le&
\sum_{j=1}^{k_2}c_j\|f_j\|+d(u^{k_1})^{-1}d(u^{k_2-k_1})^{-1}
\|a_{k_2}\|
\le \cdots \le
\sum_{j=1}^\infty c_j\|f_j\| ,
\eeq
since only a finite number of elements $a_j$ are non-zero. 
\end{proof}

Using a construction from \cite{muller} we obtain 
a similar result for extensions of Banach
space operators.

\begin{thm}\label{thm:*ban}
Let $T$ be an operator acting on a Banach space
$X$. Let $(c_j)$ be a sequence of positive numbers satisfying
condition $(*)$ for $T \in B(X)$.
Then there exists a  Banach space $Y$ containing $X$ as
a closed subspace and an invertible operator $S \in B(Y)$ such
that $S|X = T$ and
$\|S^{-j}\|\le c_j\quad(j\ge 1)$.
Moreover, we have $\|S^j\| \leq \|T^j\|\quad(j\ge 1)$ and 
$\sigma(S) \subset \sigma(T)$.
\end{thm}

\begin{proof}
Let $\AA$ be a maximal commutative subalgebra
of $B(X)$ containing $T$.  

Set $\BB = \AA \oplus X$. Define the norm and multiplication in
$\BB$ by 
$$\|A \oplus x\| = \|A\| +
\|x\|$$ 
and 
$$(A \oplus x) (A'\oplus x') =  
AA' \oplus (Ax' + A' x) \qquad (A,A' \in \AA, x,x' \in X) .$$ 
Then $\BB$ is a commutative Banach algebra and $A \mapsto A
\oplus 0 \ \ (A \in \AA)$ is an isometrical embedding $\AA \to
\BB$. 

Let $n\ge 0$. It is easy to show that
$$
d^{\BB} \big(T^n \oplus 0) =  d^{B(X)} (T^n) = m(T^n).
$$
By Theorem \ref{thm:arens}, 
there exists a commutative Banach algebra $\CCC \supset
\BB$ such that $T\oplus 0$ is invertible in $\CCC$ and
$$
\|(T \oplus 0)^{-j}\|_\CCC\le c_j\qquad(j\ge 1).
$$
Consider the operator $S:\CCC \to \CCC$ defined by  $Sc = (T
\oplus 0) c \quad  (c \in \CCC)$. Then $S$ is invertible and 
$$
\|S^{-j}\|\le c_j\qquad(j\ge 1).
$$
For $x \in X$ we have
$$
S(0 \oplus x) = (T \oplus 0) (0 \oplus x) = 0 \oplus T x.
$$
If we identify
$x \in X$ with $0 \oplus x \in \BB \subset \CCC$, then
$T=S_{\mid X}$. 

The relation $\|S^j\| \leq \|T^j\|\quad(j\ge 1)$ is easy to verify.

Finally, we have 
$$\sigma^{B(X)} (T)=\sigma^A (T)
\supset\sigma^B (T\oplus 0)
\supset \sigma^C (T\oplus 0)
\supset \sigma^{B(C)} (S) .$$
\end{proof}

\section{Extensions to $SQ_p(X)$-spaces}
In this section we study the similarity to restrictions of 
invertible operators  
acting on $SQ_p(X)$-spaces.

The proof of the following result uses an idea from \cite{batty/yeates}.
\begin{thm}\label{equivSQp}
Let $(c_j)_{j=1}^\infty$ be a sequence of positive numbers 
which is submultiplicative. Let $p\ge 1$ be a fixed real number, $X$ a Banach
space and $T\in\Bx$. 

(1) \quad  Suppose that there exists a  Banach space 
$Y$, $M\ge 1$, an operator $\pi : X \to Y$ such that 
$\|x\|\le M\|\pi(x)\|$
for all $x\in X$, 
and an invertible operator $S \in B(Y)$ such
that $S\pi = \pi T$ and $S^{-1}$ is $(c,p)$-near the null operator modulo
$\pi(X)$, that is
$$\left\| \sum_{j=1}^nS^{-j}\pi(y_j)\right\| \leq \left(
\sum_{j=1}^nc_j^p\left\|y_j\right\|^p\right)^{1/p}$$
for every $n\ge 1$ and all $y_j \in X$.
 Then we have 
$$\|x\|^p \le M^p\left(c_n^p\|x_0\|^p + c_{n-1}^p\|x_1\|^p 
+ \cdots + c_1^p\|x_{n-1}\|^p\right),$$
whenever $T^nx = x_0+Tx_1+\cdots +T^{n-1}x_{n-1}$.

(2) \quad Let $M \ge 1$ and $p\ge 1$. Suppose that the equality 
$$T^nx = x_0+Tx_1+\cdots +T^{n-1}x_{n-1} \quad  (x_i \in X, 1\le i\le n)$$ 
always implies
$$\|x\|^p \le M^p\left(c_n^p\|x_0\|^p + c_{n-1}^p\|x_1\|^p 
+ \cdots +c_1^p\|x_{n-1}\|^p\right) .$$
Then there exists a  Banach space $(Y,|\cdot|)$ which is a $SQ_p(X)$-space,  
an isomorphic embedding $\pi : X \to Y$ satisfying 
$\frac{\|x\|}{M 2^{(p-1)/p}}\le |\pi(x)|\le\|x\|\quad(x\in X)$, 
and an invertible operator $S \in B(Y)$ such
that $S\pi = \pi T$ and
$\|S^{-j}\|\le c_j$ for every $j\ge 1$. Moreover, 
$S^{-1}$ is $(c,p)$-near the null operator modulo
$\pi(X)$,
$\|S^j\| \leq \|T^j\|\quad(j\ge1)$ and $\sigma(S) \subset
\sigma(T)$. 
\end{thm}

\begin{proof}
(1) \quad Suppose that $T$ has an invertible extension 
$S$ as in the statement of the 
theorem and let $\pi : X \to Y$ satisfy $\|x\|\le M\|\pi(x)\|$
for all $x\in X$ and  
$S\pi = \pi T$. Suppose that $T^nx = x_0+Tx_1+\cdots +T^{n-1}x_{n-1}$. Then
\beq 
\|x\| &\le&M\|\pi(x)\| = M\|S^{-n}S^n\pi(x)\| = M\|S^{-n}\pi(T^nx)\|\\
&=& M\Bigl\|S^{-n}\pi\Bigl(\sum_{k=0}^{n-1}T^kx_k\Bigr)\Bigr\| 
= M\Bigl\|\sum_{k=0}^{n-1}S^{-(n-k)}\pi(x_k)\Bigr\|\\
&\le& 
M\left(\sum_{k=0}^{n-1}c_{n-k}^p\|x_k\|^p\right)^{1/p}.
\eeq

(2) \quad Suppose now that 
$$\|x\|^p \le M^p\left(c_n^p\|x_0\|^p + c_{n-1}^p\|x_1\|^p + 
\cdots + c_1^p\|x_{n-1}\|^p\right) ,$$
whenever $T^nx = x_0+Tx_1+\cdots +T^{n-1}x_{n-1}$. For $x_0=T^nx$ we get
$$\|T^nx\| \ge \frac{1}{Mc_n}\|x\| .$$
In particular, each operator $T^n$ is injective.

{\bf The equivalence relation.} Let $X_0 = X\times\Z$. 
We define an equivalence relation on $X_0$ by 
$(x,t)\sim (y,s)$ if there exists $m\in\N$ such that $s+m\in\N$, $t+m\in\N$ 
and $T^{s+m}x = T^{t+m}y$. 

Let $X_1 = X_0/\sim$ be the space of equivalence classes. We denote the
equivalence class containing $(x,t)$ by $[x,t]$. Each equivalence class contains
a member $(x,t)$ with $t\in\N$. 

{\bf The operations.} The operations 
\begin{center}
$[x,t]+[y,s] = [T^sx+T^ty,s+t]$, \quad $s,t\in \N$, $x,y\in X$,\\
$\alpha[x,t] = [\alpha x,t]$, \quad $t\in \N , \alpha \in \C$,
\end{center}
endow $X_1$ with a structure of vector space. 

{\bf The norm.} Set $c_0=1$. We define the norm on $X_1$ as
follows. For $[x,t]\in X_1$, set 
$$\left|[x,t]\right|^p = \inf \left\{ \sum_{i=0}^n \|x_i\|^pc_{i}^p : 
n\in \N,
\sum_{i=0}^n[x_i,i] = [x,t]\right\}.$$ 
We note that 
the existence of a decomposition 
$[x,t] = \sum_{i=0}^n[x_i,i]$ with $t \ge n$ is
equivalent to 
$$x = \sum_{i=0}^n T^{t-i}x_i .$$

It is easy to see that $|\cdot|$ is well-defined and 
$|\lambda[x,t]| = |\lambda||[x,t]|$ \quad 
($\lambda \in \C$).

Let $[x,t]$ and $[y,s]$ be two elements of $X_1$ decomposed by 
$[x,t]=\sum_i[x_i,i]$ and $[y,s]=\sum_i[y_i,i]$. Then
$[x,t]+[y,s]= \sum_i[x_i+y_i,i]$. By the triangular inequality
in $\ell^p$, we have
\beq
|[x,t]+[y,s]| &\le &
(\sum_i\|x_i+y_i\|^p c_i^p)^{1/p}\le
(\sum_i(\|x_i\|+\|y_i\|)^p c_i^p)^{1/p}\\
 &\le &
(\sum_i\|x_i\|^p c_i^p)^{1/p}+ (\sum_i\|y_i\|^p c_i^p)^{1/p}.
\eeq
Taking the infimum on the right hand side over all
decompositions of $[x,t]$ and $[y,s]$ we get
$|[x,t]+[y,s]|\le |[x,t]|+|[y,s]|$.

We show that $|\cdot|$ is a norm. Let $x\in X$ and $t\ge 0$. 
Consider a decomposition
$$[x,t] = \sum_{i=0}^n[x_i,i]$$
with $x_i\in X$. Then
$$[x,t] = \sum_{i=0}^n[x_i,i] = \sum_{i=0}^n[T^{n-i}x_i,n]
= \left[\sum_{i=0}^nT^{n-i}x_i,n\right].$$
Hence 
$$T^n(x-T^tx_0) = T^t\left(\sum_{i=1}^nT^{n-i}x_i\right) =
\sum_{i=1}^nT^{n-i}(T^tx_i).$$
By hypothesis, we have 
$$\|x-T^tx_0\|^p \le M^p\left(\sum_{i=1}^nc_{i}^p\|T^tx_i\|^p\right) .$$
Since 
$$\frac{1}{2^{p-1}}\|x\|^p - \|T^tx_0\|^p \le \|x-T^tx_0\|^p,$$
we get 
$$\frac{1}{2^{p-1}}\|x\|^p \le M^p\sum_{i=0}^nc_{i}^p\|T^tx_i\|^p \le 
M^p\|T^t\|^p\sum_{i=0}^nc_{i}^p\|x_i\|^p .$$
Since this is true for all such decompositions, 
we obtain 
$$|[x,t]| \ge \frac{1}{2^{(p-1)/p}M\|T^t\|}\|x\|.$$
In particular, $|[x,t]|\ne 0$ whenever $x\ne 0$.

{\bf The isomorphic embedding $\pi$.} The space $X$ embeds isomorphically into $X_1$. 
The imbedding is given by $\pi:x
\to [x,0]$ and the trivial decomposition $[x,0] = [x,0]$ gives $|\pi(x)|\le
\|x\|$. The previous paragraph, for $t=0$, shows that 
$$|\pi(x)| \ge \frac{1}{M2^{(p-1)/p}}\|x\|.$$

{\bf The operator $S$.} Define $S$ on $X_1$ by $S[x,s] =
[x,s-1]$, $x\in X$, $s\in\Z$. Clearly the definition of $S$ is
correct, $S$ is a linear map and $S\pi = \pi T$. 

The inequality 
$$\left| S^j[x,t]\right| \le \|T^j\|\cdot \left|[x,t]\right| $$
can be proved exactly as in \cite{batty/yeates}.
Thus $\|S^{j}\|\le \|T^j\|$ for all $j\ge 0$.

We show now that $|S^{-s}[x,t]| \leq c_s|[x,t]|$ for all positive $s$ and all
classes $[x,t]$. Consider a decomposition  
$$[x,t] = \sum_{i=0}^n[x_i,i]$$
with $x_i\in X$. Then 
$$[x,t+s] =  \sum_{i=0}^n[x_i,i+s].$$
Thus 
$$\left|[x,t+s]\right|^p \leq \sum_{i=0}^nc_{i+s}^p\|x_i\|^p .$$
Using the submultiplicativity of the sequence 
$c = (c_j)_{j=1}^{\infty}$ we obtain
$$
|S^{-s}[x,t]|^p=
\left|[x,t+s]\right|^p \leq c_s^p\sum_{i=0}^nc_{i}^p\|x_i\|^p .$$
This yields the announced estimate. 

We show now that
$$\left| \sum_{j=1}^nS^{-j}\pi(y_j)\right| \leq \left(
\sum_{j=1}^nc_j^p\left\|y_j\right\|^p\right)^{1/p}$$
for every $n\ge 1$ and all $y_j \in X$. Indeed, we have 
$$\sum_{j=1}^nS^{-j}\pi(y_j) = \sum_{j=1}^nS^{-j}[y_j,0] = 
\sum_{j=1}^n[y_j,j].$$
Therefore
$$\left|\sum_{j=1}^nS^{-j}\pi(y_j)\right|^p \leq  \sum_{j=1}^nc_j^p\|y_j\|^p .$$
In fact, the same arguments provide the stronger (if $p >1$) inequality
$$\left|\sum_{j=0}^nS^{-j}\pi(y_j)\right|^p \leq  \sum_{j=0}^nc_j^p\|y_j\|^p ,$$
for all $y_j \in X$, $j \ge 0$.

{\bf The space $Y$.} We take the Banach space $Y$ to be the completion of $X_1$ with the norm
$|\cdot|$ and extend $S$ continuously to an operator (also denoted by) 
$S$ on $Y$. 

We show now that $Y$ is an $SQ_p(X)$-space. Let $[a_{ij}]$ 
be an $n\times n$ matrix with
complex entries such that $\|a\|_{p,X} \leq 1$ (the definition of $\|a\|_{p,X}$
is recalled in the Introduction). Let $[x_j,t_j]$ be elements of
$X_1$ with decompositions 
$$[x_j,t_j] = \sum_{r=0}^{n^{(j)}}[w_{r}^{(j)}, r].$$
We have
\beq
&&\sum_i \left|\sum_j [a_{i,j}x_j,t_j]\right|^p =
\sum_i \left|\sum_j\sum_r [a_{i,j}w_r^{(j)},r]\right|^p\\
 &\le &
\sum_i\sum_r c_r^p \left\|\sum_j a_{i,j} w_r^{(j)}\right\|^p = 
\sum_r c_r^p \sum_i \left\|\sum_j a_{i,j} w_r^{(j)}\right\|^p\\
&\le &
\sum_r c_r^p \sum_j\left\|w_r^{(j)}\right\|^p=
\sum_j\sum_r c_r^p \left\|w_r^{(j)}\right\|^p.
\eeq

By taking the infimum over all possible decompositions, we get
$$
\sum_i \left|\sum_j [a_{i,j}x_j,t_j]\right|^p\le
\sum_j\left|[x_j,t_j]\right|^p.
$$
Thus $\|a\|_{p,Y}\le 1$, and so \cite{hernandez} 
$X_1$ and $Y$ are $SQ_p(X)$-spaces.

{\bf Spectrum behaviour.} Suppose that $T-\ld$ is invertible 
in $B(X)$. Define ${L}$ on $X_1$ by ${L}[x,t] =
[(T-\ld)^{-1}x,t]$. It is easy to see that the definition of $L$
is correct. We have
$$(S-\ld)[x,t] = [x,t-1]-[\ld x,t] = [(T-\ld)x,t]$$
and ${L}(S-\ld)[x,t]=(S-\ld)L[x,t] = [x,t]$. Hence $S-\ld$ is invertible
in $B(Y)$.
\end{proof}
\begin{rmks}
(a) The embedding $\pi$ becomes isometric if $M=p=1$ (for instance).
The case $M=p=1$ and $c_j=1$ for $j\ge
1$ was considered in \cite{batty/yeates}.

(b)
An alternative definition of the norm in $X_1$ is
$$\left|[x,t]\right|^p = \inf \left\{ \sum_{i=1}^n \|x_i\|^pc_{i}^p : 
n\in \N,
\sum_{i=1}^n[x_i,i] = [x,t]\right\}.$$ 
The difference is that decompositions of $[x,t]$ start now at $i=1$. 
The construction of $Y$, $S$ and $\pi: X\to Y$ remains
unchanged. The embedding $\pi$ satisfies in this case
$$
\frac{\|x\|}{M}\le |\pi(x)|\le c_1\|T\|\cdot \|x\|\qquad(x\in X).
$$
The remaining properties are without any change.

(c)
Note that $\sigma_{ap}(T)\subset\sigma_{ap}(S)$.

(d) We also note that
Theorem \ref{equivSQp} has a generalization to representations of 
semigroups (like in
\cite{batty/yeates}). 
\end{rmks}
\begin{dfn}Let $X$ be a Banach space, $T \in B(X)$, 
and let $p \ge 1$ be a fixed
real number. Let
$(c_j)_{j=1}^\infty$ be a sequence of positive numbers 
which is submultiplicative. We say
that $(c_j)$ \emph{satisfies condition} $(*)_p$ \emph{for} $T \in B(X)$ if there
exists an increasing sequence of integers $(k_n)$ such that 
$0=k_0<k_1< k_2<\cdots$ and 
$$
(*)_p \quad c_j\ge \frac{2^{(n+1)(p-1)/p}(k_{n+1}-k_n)^{(p-1)/p}}{m(T^{k_1}) m(T^{k_2-k_1})\cdots
m(T^{k_{n+1}-k_n})} \|T^{k_{n+1}-j}\|
$$
for all $n\ge 0$ and $j$ satisfying $k_n<j\le k_{n+1}$. 

We say
that $(c_j)$ \emph{satisfies condition} $(*)_{\infty}$ \emph{for} $T \in B(X)$ if there
exists an increasing sequence of integers $(k_n)$ such that 
$0=k_0<k_1< k_2<\cdots$ and 
$$
(*)_{\infty} \quad c_j\ge \frac{2^{n+1}(k_{n+1}-k_n)}{m(T^{k_1}) m(T^{k_2-k_1})\cdots
m(T^{k_{n+1}-k_n})} \|T^{k_{n+1}-j}\|
$$
for all $n\ge 0$ and $j$ satisfying $k_n<j\le k_{n+1}$.
\end{dfn}
The condition $(*)_1$ is the same as condition $(*)$ considered above for
Banach algebra elements. Clearly $(*)_p$ implies $(*)_q$ whenever 
$\infty\ge p\ge
q\ge 1$~; in particular, $(*)_{\infty}$ implies all other conditions
$(*)_{p}$, $p\ge 1$.
\begin{lem}\label{lem:aux}
Let $p \ge 1$ be a fixed
real number.
Suppose that $(c_j)$ is a sequence of positive numbers
satisfying condition $(*)_p$ for $T\in\Bx$. 
Then 
$$\|x\|^p \le c_m^p\|x_0\|^p + c_{m-1}^p\|x_1\|^p + \cdots + 
c_1^p\|x_{m-1}\|^p,$$
whenever $T^mx = x_0+Tx_1+\cdots +T^{m-1}x_{m-1}$.
\end{lem}

\begin{proof}
Suppose that $k_n < m \le k_{n+1}$ and that the conclusion of 
the lemma was proved for decompositions of the 
form $T^{k_{n+1}}x = x_0+Tx_1+\cdots +T^{k_{n+1}-1}x_{k_{n+1}-1}$. If 
$T^my = y_0+Ty_1+\cdots +T^{m-1}y_{m-1}$, then 
$$T^{k_{n+1}}y = 0+ \cdots + T^{k_{n+1}-m}y_0 + T^{k_{n+1}-m+1}y_1 + \cdots +
T^{k_{n+1}-1}y_{m-1}$$
and 
the lemma will be also proved for decompositions starting with $T^my$.

So suppose that 
$$T^{k_{n+1}}x = \sum_{j=1}^{k_{n+1}}T^{k_{n+1}-j}x_{k_{n+1}-j} .$$
Then, using the inequality $\|a-b\|^p \ge \frac{1}{2^{p-1}}\|a\|^p - \|b\|^p$,
we have
\beq
&&\left\|x_0\right\|^p  = \left\|T^{k_{n+1}}x - 
\sum_{j=1}^{k_{n+1}-1}T^{k_{n+1}-j}x_{k_{n+1}-j}\right\|^p\\
    & = & \left\|T^{k_{n+1}-k_n}\left(T^{k_n}x -
    \sum_{j=1}^{k_{n}}T^{k_n-j}x_{k_{n+1}-j}\right) -
    \sum_{j=k_{n}+1}^{k_{n+1}-1}T^{k_{n+1}-j}x_{k_{n+1}-j}\right\|^p\\
    &\ge & \frac{1}{2^{p-1}}m(T^{k_{n+1}-k_n})^p\left\|T^{k_n}x -
    \sum_{j=1}^{k_{n}}T^{k_n-j}x_{k_{n+1}-j}\right\|^p \\
    & - &
    \left\|\sum_{j=k_{n}+1}^{k_{n+1}-1}T^{k_{n+1}-j}x_{k_{n+1}-j}\right\|^p .
\eeq
Using now the inequality 
$$ \left\| \sum_{i=1}^N a_i\right\|^p \le 
 N^{p-1}\left( \sum_{i=1}^N\left\|a_i\right\|^p\right),$$
we obtain
\beq
&  &\|x_0\|^p +
(k_{n+1}-k_{n}-1)^{p-1}
\sum_{j=k_{n}+1}^{k_{n+1}-1}\left\|T^{k_{n+1}-j}\|^p\|x_{k_{n+1}-j}\right\|^p \\
&\ge & \frac{1}{2^{p-1}}m(T^{k_{n+1}-k_n})^p\left\|T^{k_n}x -
    \sum_{j=1}^{k_{n}}T^{k_n-j}x_{k_{n+1}-j}\right\|^p .
\eeq
Writing again 
$$T^{k_n}x -
    \sum_{j=1}^{k_{n}}T^{k_n-j}x_{k_{n+1}-j}$$ as 
$$T^{k_{n}-k_{n-1}}
\left(T^{k_{n-1}}x -
\sum_{j=1}^{k_{n-1}}T^{k_{n-1}-j}x_{k_{n+1}-j}\right)
- \sum_{j=k_{n-1}+1}^{k_{n}}T^{k_{n}-j}x_{k_{n+1}-j}
$$
and applying the same inequalities, we arrive after several steps at 
\beq
\|x_0\|^p +
(k_{n+1}-k_{n}-1)^{p-1}
\sum_{j=k_{n}+1}^{k_{n+1}-1}\left\|T^{k_{n+1}-j}\|^p\|x_{k_{n+1}-j}\right\|^p
&+ & \\
\sum_{r=1}^n\left(\frac{1}{2^{p-1}}\right)^{r}m(T^{k_{n+1}-k_{n}})^p 
\cdots m(T^{k_{n-r+2}-k_{n-r+1}})^p \left( k_{n-r+1}-k_{n-r}\right)^{p-1} & & \\
\times\sum_{j=k_{n-r}+1}^{k_{n-r+1}}\left\|
T^{k_{n-r+1}-j}\right\|^p\left\|x_{k_{n+1}-j}\right\|^p & &\\
\ge \left(\frac{1}{2^{p-1}}\right)^{n+1}m(T^{k_{n+1}-k_{n}})^p 
\cdots m(T^{k_{2}-k_{1}})^p m(T^{k_{1}})^p\left\|x\right\|^p . & & 
\eeq
This yields 
\beq
\|x\|^p\!\! &\le &\!\! \sum_{r=0}^n
\frac{(2^{p-1})^{n+1-r}(k_{n-r+1}-k_{n-r})^{p-1}}{m(T^{k_{n-r+1}-k_{n-r}})^p 
\cdots m(T^{k_{1}})^p}\!\sum_{j=k_{n-r}+1}^{k_{n-r+1}}\!\left\|
T^{k_{n-r+1}-j}\right\|^p\left\|x_{k_{n+1}-j}\right\|^p \\
 & \le & \sum_{r=0}^n
 \sum_{j=k_{n-r}+1}^{k_{n-r+1}}c_j^p\left\|x_{k_{n+1}-j}\right\|^p =
 \sum_{j=1}^{k_{n+1}}c_j^p\left\|x_{k_{n+1}-j}\right\|^p.
\eeq
\end{proof}

The above results imply the following generalization of Theorem 2.4.
\begin{thm}\label{thm:obt}
Let $p\ge 1$. Let $T$ be an operator acting on a Banach space
$X$. Let $(c_j)_{j=1}^\infty$ be a sequence of positive numbers satisfying
condition $(*)_p$ for $T \in B(X)$.
Then there exists a  Banach space $Y$ which is a $SQ_p(X)$-space, 
an isomorphic embedding $\pi : X \mapsto Y$ satisfying
$\frac{\|x\|}{2^{(p-1)/p}}\le \|\pi(x)\|\le\|x\|\quad(x\in X)$
and an invertible operator $S \in B(Y)$ such that $S\pi = \pi T$, 
$\|S^{-j}\|\le c_j\quad(j\ge1)$ and
$\|S^j\| \leq \|T^j\|\quad(j\ge1)$. Moreover, 
$S^{-1}$ is $(c,p)$-near the null operator modulo
$\pi(X)$ and $\sigma(S) \subset
\sigma(T)$. 
\end{thm}

\section{Applications}\label{sect:applications}
The previous extension results give a general way of constructing
invertible extensions of an operator with prescribed growth
conditions.  For an operator $T\in B(X)$ we write for short
$$v_n(T)=\max\{\|T^n\|,m(T^n)^{-1}\}\quad(n\ge 0).$$

We consider the following growth conditions for $T$:
\begin{itemize}
\item[$(P(s))$](polynomial growth condition) there are $C>0$ and $s\ge0$ 
such that $v_n(T)\le
Cn^s\quad(n\ge 1)$;

\item[$(B)$](Beurling-type condition) 
$\sum_{n=1}^\infty \frac{\log v_n(T)}{n^2}<\infty$;

\item[$(E(s))$](Exponential growth) there are $C>0$ and $0<s<1$ 
such that $v_n(T)\le
Ce^{n^s}\quad(n\ge 1)$.
\end{itemize}
Note that condition $(P(s))$ implies $(E(s'))$ (for any $s'>0$), which
implies $(B)$. Also, \cite{millers/neumann}, 
if $T$ satisfies $(B)$ and $T$ is invertible, then
$\sigma(T) = \sigma_{ap}(T) \subset \T$. If $T$ satisfies $(B)$ and $0\in
\sigma(T)$, then $\sigma_{ap}(T) = \T$ and $\sigma(T) = \{z: |z|\le 1\}$.

Other growth conditions can be also considered.

\subsection{$\et$-subscalar operators}
We denote by $\mathcal{E}(\C) = C^{\infty}(\C)$ the usual 
Fr\'echet algebra of all $C^{\infty}$-functions on $\C$ with the 
topology of uniform convergence of derivatives of all orders  
on  compact subsets of $\C$.
An operator $S \in \Bx$ is said \cite{colojoara/foias}
to be \emph{generalized scalar} (or $\mathcal{E}(\C)$-scalar) 
if there is a continuous algebra homomorphism
$\Phi : \mathcal{E}(\C) \to \Bx$ for which 
$\Phi(1) = I$ and $\Phi(z) = S$. A bounded linear operator 
is $\mathcal{E}(\C)$-\emph{subscalar} 
if it is similar to the restriction of a 
$\mathcal{E}(\C)$-scalar operator 
to one of its closed invariant subspaces. According to a result by J.
Eschmeier and M. Putinar (see \cite[Sect. 6.4]{eschmeier/putinar}), 
a Banach space operator 
$T$ is $\mathcal{E}(\C)$-subscalar if and only if 
$T$ has property $(\beta)_{\mathcal{E}}$, i.e., 
for every open set $U\subset \C$, the operator
$T_U$ on $\mathcal{E}(U,X)$ (the space of $C^{\infty}$-functions from $U$ into
$X$), defined by 
$T_U(f)(z) = (T-z)f(z)$, is injective and
has closed range.

The following statements are equivalent
(see \cite{colojoara/foias})~: 

(1) $T$ is $\mathcal{E}(\T)$-scalar, i.e., it 
has a continuous functional calculus on the Fr\'echet algebra 
$\mathcal{E}(\T) =C^{\infty}(\T)$ of smooth functions on the
unit circle $\T$; 

(2) $T$ is generalized scalar with $\sigma(T) \subset \T$ ;

(3) $T$ is invertible, and there exist constants $C > 0$ 
and $s \ge 0$ such that 
$$\|T^n\| \leq C(1+|n|)^s \qquad (n \in \Z ).$$

K.B.~Laursen and M.M.~Neumann \cite[Problem 6.1.15]{laursen/neumann} and 
M.~Didas \cite{didas} 
asked if 
$\mathcal{E}(\T)$-subscalar operators are characterized 
by the polynomial growth
condition $(P(s))$ above.  
We refer to \cite{didas,laursen/neumann,millers/neumann,mmn2,mmn3,mmn4} for
several partial results. By \cite{douglas} 
the hard implication holds for $s=0$ and $C=1$. 

Since
condition $(P(s))$ implies that $\sm_{ap}(T) \subset \T$,
it follows (\cite{muller,read2}) that $T$ has an invertible extension $S$ 
such that $\sm(S) =
\sm_{ap}(T) \subset \T$. 
By \cite{read3}, if $T$ acts on a Hilbert space, then $S$ acts
also on a Hilbert space. However, no control on the norms of
inverses is guaranteed by this method. 

The following result gives a 
complete positive answer.
\begin{thm}\label{subscalar}
(1) \quad 
An operator $T\in \Bx$ is $\mathcal{E}(\T)$-subscalar if and only if
there exist constants $C > 0$ and $s \ge 0$ such that
$$(P(s))\quad \frac{1}{Cn^{s}}\|x\| \leq \|T^nx\| \leq Cn^{s}\|x\|
\quad (x\in X ; n \in \N) .$$ 
Moreover, given $p\ge1$, there exist   
a $SQ_p(X)$-space $Y$, an invertible $\mathcal{E}(\T)$-scalar 
operator $S$ on $Y$ and a closed subspace $M\subset Y$ invariant
with respect to $S$ such that $T$ is similar to the restriction
$S_{\mid M}$. 
We also have $\sm(S) = \sm_{ap}(T)$.

For $p=1$ the operator $S$ is an extension of $T$.

(2) \quad If the Hilbert space operator $T \in B(H)$ verifies 
$$(P(s))\quad \frac{1}{Cn^{s}}\|h\| \leq \|T^nh\| \leq Cn^{s}\|h\|
\quad (h\in H ; n \in \N) ,$$
then there exists a Hilbert space $K$ and 
a $\mathcal{E}(\T)$-scalar extension $S\in B(K)$ with $\sm(S) =
\sm_{ap}(T)$.
\end{thm}
\begin{proof}
(1) \quad Suppose that $T$ is similar to an operator having 
a $\mathcal{E}(\T)$-scalar extension $S$. 
According to the above mentioned result, 
$S$ is $\mathcal{E}(\T)$-scalar if and
only if $S$ is invertible and $\|S^n\|$ is 
bounded by a constant times $(1+|n|)^{s}$, for each $n \in \Z$.
Therefore, restrictions of $\mathcal{E}(\T)$-scalar operators satisfy the 
growth condition $(P(s))$ from the theorem.
Consequently, $T$ satisfies $(P(s))$. 

Suppose now that $T$ satisfies the growth condition $(P(s))$.
Let $C>0$ and $s\ge0$ satisfy $v_n:= v_n(T)\le Cn^s\quad(n\ge 1)$.
Let $\e>0$. Then $\lim_{n\to\infty}\frac{v_n}{n^{s+\e/6}}=0$.
Choose $k_1\ge e^4$ such that $v_n\le n^{s+\e/6}$ 
for all $n\ge
k_1$. 

Let 
$$K=\max\{2k_1\|T^j\|\cdot 
m(T^{k_1})^{-1}:0\le j\le k_1\}$$
and set $c_j=K(j+1)^{6s+3+\e}$. 
Clearly $(c_j)$ is a submultiplicative
sequence.

We show that $(c_j)$ satisfies condition $(*)_{\infty}$ for $T$. 
Set $k_n=k_1^{2^{n-1}}\quad(n\ge 1)$.

For $j\le k_1$ we have 
$$
2k_1m(T^{k_1})^{-1}\cdot\|T^{k_1-j}\|\le K\le c_j.
$$
Let $n\ge 1$ and $k_n< j\le k_{n+1}$. Then $2^{n-1}\log k_1\le \log j$ and
\beq
&&2^{n+1}(k_{n+1}-k_n)m(T^{k_1})^{-1}\cdots
m(T^{k_{n+1}-k_n})^{-1}\|T^{k_{n+1}-j}\|\\
&\le& 
2^{n+1}
k_{n+1}\Bigl(k_1k_2\cdots k_{n+1}k_{n+1} \Bigr)^{s+\e/6} \\
&\le&
\left(\frac{2^2}{\log k_1}\log j\right)k_{1}^{2^n}
\bigl(k_1k_1^2\cdots k_1^{2^n}k_1^{2^n}\bigr)^{s+\e/6}
\\
&\le&
\left(\log
j\right)\left(k_{1}^{2^{n-1}}\right)^{2}\bigl(k_1^{3\cdot 2^n}\bigr)^{s+\e/6}\\
&\le& 
j\Bigl(k_1^{2^{n-1}}\Bigr)^{2+6s+\e}\le j^{6s+\e+3}\le c_j.
\eeq

Thus $(c_j)$ satisfies condition $(*)_{\infty}$. If $p\ge 1$ is fixed, then 
$(c_j)$ also satisfies condition $(*)_{p}$. 
By Theorem 3.5, there exists an
invertible operator $S$ on a $SQ_p(X)$-space $Y$
extending $T$ up to a similarity and satisfying $\|S^j\|=\|T^j\|$
and $\|S^{-j}\|\le c_j$
for all $j\ge 1$. Clearly $S$ has property $(P(6s+\e+3))$. Moreover, $S^{-1}$ is
$(c,p)$-near the null operator modulo $X$.

For $p=1$, the space $X$ is isometrically embedded into $Y$, and
so $S$ is an extension of $T$.

Since $\sigma(S)\subset\T$, we have $\sigma_{ap}(S)=\sigma(S)$.
By the spectral radius formula we have
$\sigma(T)\subset\{z:|z|\le 1\}$. By \cite{mz},
$$\min\{|z| : z\in\sigma_{ap}(T)\}=\lim_{n\to\infty} m(T^n)^{1/n}\ge 1 .$$ 
Thus
$\sigma_{ap}(T)=\sigma(T)\cap\T$. By Theorem 3.1,
$\sigma_{ap}(T)\subset\sigma(S)\subset\sigma(T)$. Hence 
$\sigma(S)=\sigma_{ap}(T)$.

(2) \quad Since $(c_j)$ satisfies condition $(*)_2$ for $T$, it
follows from Theorem \ref{thm:obt} that 
there exists a Hilbert space $K$, 
an isomorphic embedding $\pi : H \mapsto K$ and
an $\mathcal{E}(\T)$-scalar operator $S \in B(K)$ satisfying $S\pi = \pi T$. 
We can introduce a new equivalent Hilbert space norm on
$K$ such that $\pi$ becomes an isometry. Indeed, let $P$ be
the orthogonal projection onto $\pi (H)$. Define the new norm on
$K$ by 
$$|||u|||=( \|\pi^{-1}Pu\|_H^2+\|(I-P)u\|_K^2)^{1/2}\qquad(u\in
K).$$ 
We have $|||\pi(x)||| = \|x\|_H$ for all $x\in H$. Then $S$, acting on
the Hilbert space $(K,|||\cdot|||)$, is the required $\mathcal{E}(\T)$-scalar 
extension of $T$.
\end{proof}
\begin{rmk}
Let $H$ be the Hilbert space with an orthonormal basis
$(e_n)$ $\quad (n=0,1,\dots)$. It is easy to see that the Bergman shift on
$H$, given by
$$Be_n = \sqrt{\frac{n+1}{n+2}}e_{n+1},$$ 
satisfies the polynomial growth 
condition $(P(1/2))$. Therefore, the
Bergman shift has a generalized scalar extension with spectrum the unit
circle. This has to be compared to the known fact that $B$ is subnormal, with
minimal normal extension (the multiplication by the variable $z$ on
$L^2(\D,\mu)$, where $\mu$ is the Lebesgue measure in $\D$) 
having as spectrum the closed 
unit disk $\overline{\D}$. 
\end{rmk}

\begin{prob}
Let $s \ge 0$. What is the optimal value of $s' = f(s)$ such that 
every $T\in B(X)$ satisfying $(P(s))$ has an
invertible extension satisfying $(P(s'))$? 
What is the optimal value of $s' = g(s)$ such that 
every $T\in B(H)$ satisfying $(P(s))$ has an
invertible Hilbert space extension
satisfying $(P(s'))$~?
\end{prob}

The proof of
Theorem \ref{subscalar} can be modified 
to give, for fixed $\e > 0$ and $T\in B(X)$, a Banach space $Y$ and an 
extension (with an isometric embedding) 
$S \in B(Y)$ satisfying condition $(P(6s+\e))$. Indeed, 
with $k_1$ as in the proof of Theorem \ref{subscalar}, let
$$K=\max\{\|T^j\|\cdot m(T^{k_1})^{-1}:0\le j\le k_1\}$$
and set $c_j=K(j+1)^{6s+\e}$. 
Then a similar proof shows that the sequence $(c_j)$
satisfies condition $(*)_1$ for $k_n=k_1^{2^{n-1}}\quad(n\ge 1)$.

We also notice that $g(0) = 0$. Indeed, if a Hilbert space operator $T\in B(H)$
satisfies $(P(0))$, then by \cite{sznagy} there exists an invertible operator
$L\in B(H)$ such that $V=L^{-1}TL$ is an isometry. Let $U$ be a unitary
extension of $V$ on a larger Hilbert space $K=H\oplus H^{\perp}$. 
Then $(L\oplus I)U(L\oplus
I)^{-1}$ is an extension of $T$ satisfying $(P(0))$. 
\bigskip

We can consider representations of $\N^n$ to deal with 
$\mathcal{E}(\T^n)$-subscalar operators. The proof of the following result
follows a different approach.
\begin{thm}
An $n$-tuple of commuting Banach space operators is
$\mathcal{E}(\T^n)$-subscalar if and only if each of the $n$
operators is $\mathcal{E}(\T)$-subscalar.
\end{thm}
\begin{proof}
The previous characterization of $\mathcal{E}(\T)$-subscalar
operators implies that if $T_1, \ldots ,T_n$ are commuting
$\mathcal{E}(\T)$-subscalar operators, then the product operator $T_1\cdots
T_n$ is also $\mathcal{E}(\T)$-subscalar. The result follows from
\cite[Th. 2.2.7]{didas}.
\end{proof}

\subsection{Operators with Bishop's property $(\beta)$}
Recall  
that an equivalent  definition of 
decomposable operators is the following~: 
$T\in\Bx$ is \emph{decomposable} if for every open cover 
$\C = U\cup V$, there are 
closed invariant (for $T$) subspaces $Y$ and $Z$ of $X$ such that $X = Y+Z$ and
$\sigma(T\mid Y)\subset U$, $\sigma(T\mid Z)\subset V$. 
We refer for instance to~\cite{colojoara/foias} and~\cite{laursen/neumann}.
An operator $T \in \Bx$ has
\emph{Bishop's property} $(\beta)$ if, for every open set $U\subset \C$, the operator
$T_U$ defined by $T_U(f)(z) = (T-z)f(z)$ 
on the set $\mathcal{O}(U,X)$ of holomorphic functions from $U$ into $X$ 
is injective and
has closed range.
According to a result by E. Albrecht and J. Eschmeier (see
\cite{laursen/neumann,eschmeier/putinar}), 
$T \in \Bx$ is \emph{subdecomposable} (i.e., $T$ is similar to the restriction of a
decomposable operator) if and only if $T$ has
Bishop's property $(\beta)$. 

It was proved in \cite[5.3.2]{colojoara/foias} that an invertible operator 
$S\in B(X)$ is decomposable provided that 
$$\sum_{n=-\infty}^{\infty}\frac{\log \|S^n\|}{1+n^2} < \infty .$$

The following result answers in the affirmative a question from
\cite{millers/neumann}.
 
\begin{thm}\label{beta}
Let $T \in \Bx$ be a Banach space operator such that
$$ \sum_{n=1}^{\infty}
\frac{\log \max(\|T^n\|,m(T^n)^{-1})}{n^2} < \infty .$$ 
Then there exists a Banach space
$Y\supset X$ and an invertible operator $S\in B(Y)$ such that
$T=S_{\mid X}$ and $S$ satisfies 
$$\sum_{n={-\infty}}^\infty \frac{\log \|S^n\|}{1+n^2}<\infty .$$
In particular, $T$
has Bishop's property $(\beta)$. 
Moreover, $\sigma(S)=\sigma_{ap}(T)=\sigma(T)\cap\T$.

If $X=H$ is a Hilbert space, 
then $Y=K$ can be chosen to be a Hilbert space too.
\end{thm}

\begin{proof}
Let $T\in B(X)$ satisfy (B). 
By Theorem \ref{thm:obt}, it is sufficient to show the existence of a submultiplicative sequence 
$(d_n)$
satisfying 
$\sum_{n=1}^\infty \frac{\log d_n}{n^2}<\infty$
and the condition $(*)_{\infty}$ for $T$.

Write $r_n=v_{2^n}\quad(n\ge 0).$ Clearly $r_{n+1}\le r_n^2$ for
all $n$.

\medskip
\noindent{\bf Claim a.}
$\sum_{n=0}^\infty\frac{\log r_n}{2^n}<\infty$.
\smallskip

\noindent{\bf Proof.}
Fix $n\ge 2$. For $1\le j\le 2^{n-3}$ we have
$$
v_{2^n}\le v_{2^{n-1}+j}\cdot v_{2^{n-1}-j}.
$$
Thus $\log r_{n}\le\log v_{2^{n-1}+j}+\log v_{2^{n-1}-j}$ and
$$
\frac{\log r_n}{2^{2n}}\le
\frac{\log v_{2^{n-1}+j}}{2^{2n}}+\frac{\log v_{2^{n-1}-j}}{2^{2n}}\le 
\frac{\log v_{2^{n-1}+j}}{(2^{n-1}+j)^2}
+\frac{\log v_{2^{n-1}-j}}{(2^{n-1}-j)^2}.
$$
Hence
$$
2^{n-3}\cdot \frac{\log r_n}{2^{2n}}\le
\sum_{j=1}^{2^{n-3}}\Bigl(\frac{\log v_{2^{n-1}+j}}{(2^{n-1}+j)^2}
+\frac{\log v_{2^{n-1}-j}}{(2^{n-1}-j)^2}\Bigr)
$$
and
$$
\frac{1}{8}\sum_{n=2}^\infty \frac{\log r_n}{2^n}\le
\sum_{j=1}^\infty\frac{\log v_j}{j^2}<\infty.
$$
\bigskip

Let $n$ be a non-negative integer and let 
$n=\sum_{j=0}^\infty \al_j2^j$, where $\al_j\in \{0,1\}$,  
be its binary representation.
Define 
$$b_n=\prod_{j=0}^\infty r_j^{\al_j} \mbox{ , }
c_n=\max\{b_j^2:n\le j\le 2n\} \mbox{ and } d_n = 4n^2c_n .$$
\medskip

\noindent{\bf Claim b.}
$(b_n)$ is submultiplicative, i.e., $b_{n+m}\le b_nb_m$ for all
$m,n\ge 0$.
\smallskip

\noindent{\bf Proof.}
Let $n=\sum_{j=0}^\infty\al_j2^j$ and
$m=\sum_{j=0}^\infty\beta_j2^j$ be the binary representations of
$n$ and $m$, respectively.

By induction on $j_0$, we prove the following statement:

\noindent
\begin{quote}
There are numbers $\gamma_j \quad(0\le j)$ such that
$n+m=\sum_{j=0}^\infty \gamma_j2^j$, $\gamma_j\in\{0,1\}\quad(j<
j_0)$, $\gamma_{j_0}\in\{0,1,2,3\}$,
$\gamma_j\in\{0,1,2\}\quad(j> j_0)$ and
$b_nb_m\ge\prod_{j=0}^\infty r_j^{\gamma_j}$.
\end{quote}

For $j_0=0$ the statement is clear for the numbers
$\gamma_j=\al_j+\beta_j$. 

Suppose that the statement is true for some $j_0$. We show it
for $j_0+1$.
If $\gamma_{j_0}\le 1$ then the statement is clear. Let
$\gamma_{j_0}\in\{2,3\}$. Then 
$$
n+m=\sum_{j=0}^\infty\gamma_j'2^j ,
$$
where
$\gamma'_j=\gamma_j\quad(j\ne j_0,j_0+1)$,
$\gamma'_{j_0}=\gamma_{j_0}-2$ and
$\gamma'_{j_0+1}=\gamma_{j_0+1}+1$. 
Then
$$
b_nb_m\ge\prod_{j=0}^\infty r_j^{\gamma_j}\ge
\prod_{j=0}^\infty r_j^{\gamma'_j}.
$$
The statement for $j_0> \log_2(n+m)$ gives the
inequality $b_nb_m\ge b_{n+m}$.
\medskip

\noindent{\bf Claim c.}
$(d_n)$ is submultiplicative.
\smallskip

\noindent{\bf Proof.}
Notice that $16m^2n^2 \ge 4(m+n)^2$ for all positive integers $m$ and $n$. 
We have 
\beq
d_nd_m&\ge&
4(m+n)^2\max\{b_i^2b_j^2: n\le i\le 2n, m\le j\le 2m\}\\
&\ge&
4(m+n)^2\max\{b_l^2: n+m\le l\le 2(n+m)\}=d_{n+m}.
\eeq
\medskip

\noindent{\bf Claim d.}
$\sum_{n=1}^\infty \frac{\log d_n}{n^2}<\infty$.
\smallskip

\noindent{\bf Proof.}
It is sufficient to show the analogue claim for the sequence $(c_n)$. 

For $2^j\le n<2^{j+1}$ we have
$c_n=b_i^2$ for some $i$, $i\le 2n< 2^{j+2}$. So
$c_n\le b^2_{2^{j+2}-1}=\prod_{i=0}^{j+1} r_i^2$. Thus
$$
\sum_{n=2}^\infty \frac{\log c_n}{n^2}\le \sum_{j=1}^\infty
\frac{2^j\sum_{i=0}^{j+1} 2\log r_i}{2^{2j}}\le
2\sum_{i=0}^\infty \log r_i\cdot \sum_{j=i-1}^\infty 2^{-j}\le
8\sum_{i=0}^\infty\frac{\log r_i}{2^i}<\infty.
$$
\medskip

\noindent{\bf Claim e.}
$(d_n)$ satisfies condition $(*)_{\infty}$ for $T$.
\smallskip

\noindent{\bf Proof.}
Set $k_n=2^n-1$. For $k_n<j\le k_{n+1}$ we have $2^n\le
j<2^{n+1}$, and so
$c_j\ge\prod_{i=0}^{n}r_i^2$. Hence
\beq
&&2^{n+1}(k_{n+1}-k_n)\Bigl(m(T^{k_1})m(T^{k_2-k_1})\cdots
m(T^{k_{n+1}-k_n})\Bigr)^{-1}\|T^{k_{n+1}-j}\|\\
&\le& 
2(2^n)^2r_0r_1\cdots r_n\cdot b_{k_{n+1}-j}\\
 &\le&
2j^2\prod_{i=0}^{n} r_i^2\\
&\le& d_j.
\eeq
The inequality for
$j=1=k_1$ is clear.
\medskip

Thus $(d_n)$ also satisfies condition $(*)_1$, and so there is
an invertible extension $S$ of $T$ such that $\|S^{-n}\|\le
d_n\quad(n>0)$. Hence $S$ is decomposable.

The equalities $\sigma(S)=\sigma_{ap}(T)=\sigma(T)\cap\T$ can
be shown as in Theorem 4.1.

If $X=H$ is a Hilbert space, then the sequence $(d_n)$ satisfies
condition $(*)_2$. By Theorem 3.5, there is a Hilbert space $K$,
an invertible operator $S\in B(K)$ and an isomorphic embedding
$\pi:H\to K$ with $\pi T=S\pi$ and $\|S^{-n}\|\le d_n\quad(n>0)$.
As in the proof of Theorem 4.1, $K$ can be given a new
equivalent hilbertian norm such that $\pi$ becomes an isometry.
\end{proof}

\subsection{Condition $(E(s))$}
The following consequence of Theorem \ref{beta} implies
that condition (b) from 
\cite[Th. 3.2]{mmn3} is superfluous.
\begin{cor}
Let $T \in \Bx$ satisfying the exponential condition $(E(s))$, that is, 
there are $C>0$ and $0<s<1$ such that $v_n(T)\le
Ce^{n^s}\quad(n\ge 0)$. Then $T$ has property $(\beta)$.
\end{cor}

The following result 
answers an open question from \cite{millers/neumann}.
\begin{thm}\label{thm:4.7}
Let $T\in B(X)$ satisfy (E). Then there exist a Banach 
space $Y\supset X$ 
and an invertible operator $S$ on a larger space such that
$T$ is a restriction of $S$ and $S$ satisfies 
$(E(s'))$ for suitable $s'<1$. The construction is hilbertian.
\end{thm}

\begin{proof}
Let $\e$ be an arbitrary positive number. Set $k_n=2^n\quad(n\ge
1)$. It is now a matter of routine to verify that the sequence
$c_j=K\cdot e^{j^{s+\e}}$ satisfies condition $(*)_{\infty}$ for $T$, where
$K$ is a suitable constant. Thus $T$ can be extended to an
invertible operator satisfying condition $(E(s+\e))$. 
The construction is hilbertian in the sense that
if $X=H$ is Hilbert, then $Y=K$ can be chosen a Hilbert space too.
We omit the details.
\end{proof}

\subsection{A hilbertian counterpart of Arens' result.} 
We obtain the following hilbertian counterpart of Arens' result.
\begin{cor}
Let $T\in B(H)$ be an operator on Hilbert space with $m(T) > 0$. Then there
exist a Hilbert space $K$, 
an isometric embedding $\pi :H \mapsto K$
and an invertible operator $S\in
B(K)$ such that $S\pi = \pi T$, 
$\|S^j\| \leq \|T^j\|\quad (j\ge 1)$, $\|S^{-1}\| \le 2/m(T)$
and 
$$\left\| \sum_{j=0}^NS^{-j}\pi(x_j)\right\|^2 \leq 
2\sum_{j=0}^N\left(\frac{\sqrt{2}}{m(T)}\right)^{2j}\|x_j\|^2
$$
for every $N\in \N$ and all $x_j \in H$.
\end{cor}

\begin{proof}
Let $c_j = \left(\frac{\sqrt{2}}{m(T)}\right)^j$, $j\ge 1$. Then the sequence
$(c_j)$ satisfies the condition $(*)_2$ for $T$ (take $k_n=n$). 
It follows
from the proof of Theorem \ref{thm:obt} that there
exist a Hilbert space $K$, 
an isomorphic embedding $\pi :H \mapsto K$ 
satisfying $\frac{1}{\sqrt{2}}\|x\| \le \|\pi(x)\| \le \|x\|$ 
for any $x\in H$, 
and an invertible operator $S\in
B(K)$ such that $S\pi = \pi T$, 
$\|S^j\| \leq \|T^j\|\quad (j\ge 1)$, $\|S^{-1}\| \le \sqrt{2}/m(T)$
and 
$$\left\| \sum_{j=0}^NS^{-j}\pi(x_j)\right\|^2 \leq 
\sum_{j=0}^N\left(\frac{\sqrt{2}}{m(T)}\right)^{2j}\|x_j\|^2
$$
for every $N\in \N$ and all $x_j \in H$. 
 We now introduce a new equivalent Hilbert space norm on
$K$ such that $\pi$ becomes an isometry as in the 
proof of Theorem~\ref{subscalar}. So let $P$ be
the orthogonal projection onto $\pi H$ and define the new norm on
$K$ by 
$$|||x|||=( \|\pi^{-1}Px\|_H^2+\|(I-P)x\|_K^2)^{1/2}.$$ 
Then $|||x|||^2 \le 2\|Px\|_K^2+\|(I-P)x\|_K^2 \le 2\|x\|_K^2$. In the same way
a lower bound can be obtained; we get
$\|x\| \le |||x||| \le \sqrt{2}\|x\|$ for every $x\in K$. 
Then $S$, acting on
the Hilbert space $(K,|||\cdot|||)$, verifies the required 
inequalities.
\end{proof}

\subsection{Operators with countable spectrum.} 
In the following two results we assume that the spectrum of $T$ is countable.
We refer to \cite{batty/yeates,BBG,KN} and their references for related results.
\begin{thm}
Let $T\in \Bx$ be a Banach space operator. Suppose that there are 
positive constants 
$M>0$, $C>0$ and $0< s < \frac{1}{2}$ such that 
$$\frac{1}{Ce^{n^s}}\|x\| \le \|T^nx\| \le M\|x\|$$
for every $x \in X$ and $n\in\N$. 
Suppose also that the spectrum $\sigma(T)$ of $T$ is countable. Then 
$$
(P(0)) \qquad \frac{1}{M}\|x\| \le \|T^nx\| \le M\|x\|
$$
for every $x \in X$. In particular, $T$ is $\mathcal{E}(\T)$-subscalar.
\end{thm}
\begin{proof}
We have $\|T^n\| \le M$ and $m(T^n))^{-1} \leq Ce^{n^s}$. Let $\e > 0$ be 
a positive number such
that $s+\e <\frac{1}{2}$. 
Using (the proofs of) Theorems \ref{thm:4.7} and \ref{thm:*ban} (or
\ref{thm:obt}), 
there exists a constant $K
>0$ such that $T$ has an invertible
extension $S$ on a Banach space $Y$ verifying $\|S^n\| \le M$ and
$\|S^{-n}\| \le K\exp(n^{s+\e})$ for all $n \in \N$. Moreover, 
it is possible to have an extension satisfying 
$\sigma(S) \subset \sigma(T)$. We obtain in particular that 
$$\lim_{n\to\infty}\frac{\log \|S^{-n}\|}{\sqrt{n}} = 0$$
and that the spectrum $\sigma(S)$ of $S$
is countable. From \cite[Remarque 2, p. 259]{zarrabi} we obtain 
$\|S^p\| \leq M$ for all $p\in \Z$. This yields $m(T^n)^{-1} \leq M$ for $n\ge
1$ and the
stated inequality $(P(0))$.
\end{proof}

We obtain the following
consequence in the case of Hilbert space operators. 
\begin{cor}
Let $T\in \Bh$ be a power bounded operator on a Hilbert space $H$. 
Suppose that there are positive
constants $C$ and $s < \frac{1}{2}$ such that 
$$m(T^n)^{-1} \le Ce^{n^s} \quad (n\ge 1)$$
and that $\sigma(T)$ is countable. Then $T$ is similar to a unitary 
operator.
\end{cor}
\begin{proof}
By the previous theorem, the 
operator $T$ satisfies $(P(0))$ on $H$, a condition 
which characterizes
Hilbert space operators similar to isometries \cite{sznagy}. As the spectrum is
a similarity invariant, $T$ is similar to an isometry with a countable
spectrum. Since the spectrum of a non-invertible isometry is 
the entire closed unit disk, we
obtain that $T$ is similar to a unitary operator.
\end{proof}
\subsection{Contractions with spectrum a Carleson set.}
Recall that a closed set $E$ of $\T$ is said to be a \emph{Carleson set} if 
$$ \int_{0}^{2\pi} 
\log \left( \frac{2}{\mbox{ dist }(e^{it},E)}\right)\, dt < +\infty.$$ 
\begin{thm}
Let $T\in B(H)$ be a Hilbert space contraction such that 
$\sm_{ap}(T) \subset
\T$ is a Carleson set. Suppose that there exist $C >0$ and $s\ge 0$ 
such that $m(T^n)^{-1} \le Cn^s$. Then $T$ is an isometry. 
\end{thm}
\begin{proof}
Using Theorem \ref{subscalar}, (2), 
there exist $K>0$, $s'\ge0$, a 
Hilbert space $K$ and an invertible operator $S\in B(K)$ 
which is an extension of $T$ such that
$\|S\| \le 1$, $\|S^{-n}\| \le Kn^{s'}$ and $\sm(S) = \sm_{ap}(T)$.
We obtain in particular that $\sm(S) = \sm_{ap}(T)$ is a
Carleson set. By a theorem of Esterle 
\cite{esterle} (see also \cite{kellay}), $S$ is unitary. Therefore its
restriction $T$ is an isometry.
\end{proof} 

Several results for unitaries (or operators similar to unitaries) can be
transferred to results for isometries (or operators similar to isometries) 
in an analogous manner.


\begin{thebibliography}{MMN9}
\bibitem[Ar]{arens} Arens, R. Inverse-producing extensions of normed algebras,
{\sl Trans. Amer. Math. Soc.} 88(1958), 536-548. 
\bibitem[Ba1]{badea} Badea, C. 
Perturbations of operators similar to contractions and the commutator 
equation, 
{\sl Studia Math.} 150(2002), 273--293.
\bibitem[Ba2]{badea2}Badea, C.
Operators near completely polynomially dominated ones and similarity problems, 
{\sl J. Operator Theory} 49(2003), 3--23.
\bibitem[BBG]{BBG} Batty, C. J. K.; Brze\'zniak, Z.; Greenfield, D.A.
A quantitative asymptotic theorem for contraction semigroups with countable
unitary spectrum,
{\sl Studia Math.} 121(1996), 167--183.
\bibitem[BY]{batty/yeates} Batty, C. J. K.; Yeates, S. B.
Extensions of semigroups of operators,
{\sl J. Operator Theory} 46(2001), 139--157.
\bibitem[BP]{bp}
Bercovici, H.; Petrovi\'c, S.
Generalized scalar operators as dilations, 
{\sl Proc. Amer. Math. Soc.} 123(1995), 2173--2180.
\bibitem[GT]{gindler} Gindler, H. A.; Taylor, A. E. 
The minimum modulus of a linear operator and its use in spectral theory, 
{\sl Studia Math.} 22(1962/1963), 15--41.
\bibitem[CF]{colojoara/foias} Colojoar\u{a}, I.; Foia\c{s}, C.
{\sl Theory of generalized spectral operators}, 
Mathematics and its Applications, Vol. 9. Gordon
and Breach, Science Publishers, New York-London-Paris, 1968.
\bibitem[Do]{douglas} Douglas, R.G.
On extending commutative semigroups of isometries, 
{\sl Bull. London Math. Soc.} 1 (1969), 157--159.
\bibitem[Di]{didas} Didas, M.
$\mathcal{E}(\T^n)$-subscalar 
$n$-tuples and the Cesaro operator on $H^p$, 
{\sl Ann. Univ. Sarav. Ser. Math.} 10 (2000), 
no. 2, pp. i--iii and 284--335.
\bibitem[EP]{eschmeier/putinar}Eschmeier, J.; Putinar, M. 
{\sl Spectral decompositions and analytic sheaves}, 
London Mathematical Society Monographs. New
Series, 10.
The Clarendon Press, Oxford University Press, New York, 1996.
\bibitem[Es]{esterle} Esterle, J.
Uniqueness, strong forms of uniqueness and negative powers of contractions, 
in: {\sl Functional analysis and operator theory (Warsaw, 1992)}, 127--145, 
Banach Center Publ., 30, 
Polish Acad. Sci., Warsaw, 1994. 
\bibitem[He]{hernandez}
Hernandez, R. 
Espaces $L\sp{p}$, factorisation et produits tensoriels dans les 
espaces de Banach,
{\sl C. R. Acad. Sci. Paris S\'er. I Math.}  296(1983), 
385--388.
\bibitem[Ke]{kellay} Kellay, K.
Contractions et hyperdistributions \`a spectre de Carleson, 
{\sl J. London Math. Soc.} (2) 58(1998), 185--196.
\bibitem[KN]{KN} K\'erchy, L.; van Neerven, J.
Polynomially bounded operators whose spectrum on the 
unit circle has measure
zero, 
{\sl Acta Sci. Math. (Szeged)} 63(1997), 551--562.
\bibitem[Kw]{kwapien} Kwapie\'n, S. 
On operators factorizable through $L\sb{p}$ space,
Actes du Colloque d'Analyse Fonctionnelle de Bordeaux 
(Univ. de Bordeaux, 1971), pp. 215--225. 
{\sl Bull. Soc. Math. France, Mem.} No. 31--32,
Soc. Math. France, Paris, 1972. 
\bibitem[LN]{laursen/neumann} Laursen, K.B.; Neumann, M.M.
{\sl An introduction to local spectral theory}, 
London Mathematical Society Monographs.
New Series, 20. 
The Clarendon Press, Oxford University Press, New York, 2000. 
\bibitem[LM]{lemerdy} Le Merdy, C. 
Factorization of $p$-completely bounded multilinear maps, 
{\sl Pacific J. Math.} 172 (1996), 187--213.
\bibitem[MZ]{mz} E. Makai, J. Zem\'anek, The surjectivity radius,
packing numbers and boundedness below, {\sl Integral Equations
Operator Theory} 6 (1983), 372--384.
\bibitem[MMN1]{millers/neumann} Miller,T.L. ; Miller,V. ; Neumann, M.M. 
Growth conditions and decomposable extensions, 
in: {Proceedings of the Conference on Trends in 
Banach Spaces and Operator Theory}, {\sl Contemp. Math.} 321 , 197-205,
Amer.Math.Soc. 2003.
\bibitem[MMN2]{mmn2}Miller,T.L.; Miller,V.; Neumann, M.M.
Spectral subspaces of subscalar and related operators, 
{\sl Proc. Amer. Math. Soc.} 132(2004), 1483-1493.
\bibitem[MMN3]{mmn3}Miller,T.L. ; Miller,V. ; Neumann, M.M.
Local spectral properties of weighted shifts, {\sl J.
Operator Theory} 51 (2004), 71-88.
\bibitem[MMN4]{mmn4}Miller,T.L.; Miller,V.; Neumann, M.M.
Localization in the spectral theory of operators on Banach spaces,  
in: {Proceedings of the Fourth Conference on Function 
Spaces at Edwardsville}, {\sl Contemp. Math.} 328 , 247-262, 
Amer.Math.Soc. 2003..
\bibitem[M\"u]{muller} M\"uller, V.
Adjoining inverses to noncommutative Banach algebras 
and extensions of operators, 
{\sl Studia Math.} 91 (1988), 73--77.
\bibitem[Pi]{pisier:indiana}
Pisier, G. 
Completely bounded maps between sets of Banach space operators,
{\sl Indiana Univ. Math. J.} 39(1990), 249-277.
\bibitem[Re1]{read1} 
Read, C. J. 
Inverse producing extension of a Banach algebra which 
eliminates the residual spectrum of one element, 
{\sl Trans. Amer. Math. Soc.} 286 (1984), 715--725.
\bibitem[Re2]{read2} 
Read, C. J. 
Spectrum reducing extension for one operator on a Banach space, 
{\sl Trans. Amer. Math. Soc.} 
308 (1988), 413--429.
\bibitem[Re3]{read3} 
Read, C. J. 
Extending an operator from a Hilbert space to a larger Hilbert space, 
so as to reduce its spectrum, 
{\sl Israel J. Math.} 
57 (1987), 375--380.
\bibitem[St]{stroescu}
Stroescu, E.
Isometric dilations of contractions on Banach spaces, 
{\sl Pacific J. Math.} 47(1973), 257--262.
\bibitem[SN]{sznagy}
Sz.-Nagy, B. On uniformly bounded linear 
transformations in Hilbert space, 
{\sl Acta Sci. Math. Szeged} 11(1947), 152-157.
\bibitem[SNF]{szfo} 
Sz.-Nagy, B.; Foia\c{s}, C.  
{\sl Harmonic analysis of operators on Hilbert space}, 
Translated from the French and revised,
North-Holland, Amsterdam, 1970.
\bibitem[Za]{zarrabi} Zarrabi, M.
Contractions \`a spectre d\'enombrable et propri\'et\'es d'unicit\'e 
des ferm\'es d\'enombrables du cercle, 
{\sl Ann. Inst. Fourier (Grenoble)} 43 (1993), 251--263.

\end{thebibliography}
\end{document}